\documentclass[a4paper,12pt]{article}

 \usepackage{latexsym,amssymb,amsmath}

 \newtheorem{theorem}{Theorem}
 
 \newtheorem{proposition}{Proposition}
 \newtheorem{corollary}{Corollary}
 \newtheorem{definition}{Definition}

\begin{document}

\title{Fischer decomposition by inframonogenic functions\thanks{accepted for publication in CUBO, A Mathematical Journal}}

\author{Helmuth R. Malonek$^{\star,1}$, Dixan Pe\~na Pe\~na$^{\star,2}$\\and Frank Sommen$^{\dagger,3}$}

\date{\normalsize{$^\star$Department of Mathematics, Aveiro University,\\3810-193 Aveiro, Portugal\\
$^\dagger$Department of Mathematical Analysis, Ghent University,\\9000 Gent, Belgium}\\\vspace{0.4cm}
\small{$^1$e-mail: hrmalon@ua.pt\\
$^2$e-mail: dixanpena@ua.pt; dixanpena@gmail.com\\
$^3$e-mail: fs@cage.ugent.be}}

\maketitle

\begin{abstract}
\noindent Let $\partial_{\underline x}$ denote the Dirac operator in $\mathbb R^m$. In this paper, we present a refinement of the biharmonic functions and at the same time an extension of the monogenic functions by considering the equation $\partial_{\underline x}f\partial_{\underline x}=0$. The solutions of this \lq\lq sandwich" equation, which we call inframonogenic functions, are used to obtain a new Fischer decomposition for homogeneous polynomials in $\mathbb R^m$.\vspace{0.2cm}\vspace{0.1cm}
\end{abstract}

\section{Introduction}

Let $\mathbb{R}_{0,m}$ be the $2^m$-dimensional real Clifford algebra constructed over the orthonormal basis $(e_1,\ldots,e_m)$ of the Euclidean space $\mathbb R^m$ (see \cite{Cl}). The multiplication in $\mathbb{R}_{0,m}$ is determined by the relations $e_je_k+e_ke_j=-2\delta_{jk}$ and a general element of $\mathbb{R}_{0,m}$ is of the form $a=\sum_Aa_Ae_A$, $a_A\in\mathbb R$, where for $A=\{j_1,\dots,j_k\}\subset\{1,\dots,m\}$, $j_1<\cdots<j_k$, $e_A=e_{j_1}\dots e_{j_k}$. For the empty set $\emptyset$, we put $e_{\emptyset}=1$, the latter being the identity element.

Notice that any $a\in\mathbb{R}_{0,m}$ may also be written as  $a=\sum_{k=0}^m[a]_k$ where $[a]_k$ is the projection of $a$ on $\mathbb R_{0,m}^{(k)}$. Here $\mathbb R_{0,m}^{(k)}$ denotes the subspace of $k$-vectors defined by
\[\mathbb R_{0,m}^{(k)}=\biggl\lbrace a\in\mathbb R_{0,m}:\;a=\sum_{\vert A\vert=k}a_Ae_A,\quad a_A\in\mathbb R\biggr\rbrace.\]
In particular, $\mathbb R_{0,m}^{(1)}$ and $\mathbb R_{0,m}^{(0)}\oplus\mathbb R_{0,m}^{(1)}$ are called, respectively, the space of vectors and paravectors in $\mathbb R_{0,m}$. Observe that $\mathbb R^{m+1}$ may be naturally identified with $\mathbb R_{0,m}^{(0)}\oplus\mathbb R_{0,m}^{(1)}$ by associating to any element $(x_0,x_1,\ldots,x_m)\in\mathbb R^{m+1}$ the paravector $x=x_0+\underline x=x_0+\sum_{j=1}^mx_je_j$.

Conjugation in $\mathbb R_{0,m}$ is given by
\[\overline a=\sum_Aa_A\overline e_A,\quad\overline e_A=(-1)^{\frac{\vert A\vert(\vert A\vert+1)}{2}}e_A.\]
One easily checks that $\overline{ab}=\overline b\overline a$ for any $a,b\in\mathbb R_{0,m}$. Moreover, by means of the conjugation a norm $\vert a\vert$ may be defined for each $a\in\mathbb R_{0,m}$ by putting
\[\vert a\vert^2=[a\overline a]_0=\sum_Aa_A^2.\]
The $\mathbb R_{0,m}$-valued solutions $f(\underline x)$ of $\partial_{\underline x}f(\underline x)=0$, with $\partial_{\underline x}=\sum_{j=1}^me_j\partial_{x_j}$ being the Dirac operator, are called left monogenic functions (see \cite{BDS,DSS}). The same name is used for null-solutions of the operator $\partial_x=\partial_{x_0}+\partial_{\underline x}$ which is also called generalized Cauchy-Riemann operator.

In view of the non-commutativity of $\mathbb R_{0,m}$ a notion of right monogenicity may be defined in a similar way by letting act the Dirac operator or the generalized Cauchy-Riemann operator from the right. Functions that are both left and right monogenic are called two-sided monogenic.

One can also consider the null-solutions of $\partial_{\underline x}^k$ and $\partial_x^k$ ($k\in\mathbb N$) which gives rise to the so-called $k$-monogenic functions (see e.g. \cite{F1,F2,R}).

It is worth pointing out that $\partial_{\underline x}$ and $\partial_x$ factorize the Laplace operator in the sense that
\[\Delta_{\underline x}=\sum_{j=1}^m\partial_{x_j}^2=-\partial_{\underline x}^2,\quad\Delta_x=\partial_{x_0}^2+\Delta_{\underline x}=\partial_x\overline\partial_x=\overline\partial_x\partial_x.\]
Let us now introduce the main object of this paper.

\begin{definition}
Let $\Omega$ be an open set of $\,\mathbb R^m$ $($resp. $\mathbb R^{m+1}$$)$. An $\mathbb R_{0,m}$-valued function $f\in\mathcal{C}^2(\Omega)$ will be called an inframonogenic function in $\Omega$ if and only if it fulfills in $\Omega$ the \lq\lq sandwich" equation
\begin{equation*}
\partial_{\underline x}f\partial_{\underline x}=0\quad(\text{resp.}\;\partial_xf\partial_x=0).
\end{equation*}
\end{definition}
Here we list some motivations for studying these functions.

\begin{enumerate}
\item If a function $f$ is inframonogenic in $\Omega\subset\mathbb R^m$ and takes values in $\mathbb R$, then $f$ is harmonic in $\Omega$.
\item The left and right monogenic functions are also inframonogenic.
\item If a function $f$ is inframonogenic in $\Omega\subset\mathbb R^m$, then it satisfies in $\Omega$ the overdetermined system $\partial_{\underline x}^3f=0=f\partial_{\underline x}^3$. In other words, $f$ is a two-sided 3-monogenic function.
\item Every inframonogenic function $f\in\mathcal{C}^4(\Omega)$ is biharmonic, i.e. it satisfies in $\Omega$ the equation $\Delta_x^2f=0$ (see e.g. \cite{BG,GK,M,So}).
\end{enumerate}
The aim of this paper is to present some simple facts about the inframonogenic functions (Section 2) and establish a Fischer decomposition in this setting (Section 3).

\section{Inframonogenic functions: simple facts}

It is clear that the product of two inframonogenic functions is in general not inframonogenic, even if one of the factors is a constant.

\begin{proposition}
Assume that $f$ is an inframonogenic function in $\Omega\subset\mathbb R^m$ such that $e_jf$ $($resp. $fe_j$$)$ is also inframonogenic in $\Omega$ for each $j=1,\dots,m$. Then f is of the form
\[f(\underline x)=c\underline x+M(\underline x),\]
where $c$ is a constant and $M$ a right $($resp. left$)$ monogenic function in $\Omega$.
\end{proposition}
\textit{Proof.} The proposition easily follows from the equalities
\[\partial_{\underline x}\big(e_jf(\underline x)\big)\partial_{\underline x}=-2\partial_{x_j}f(\underline x)\partial_{\underline x}-e_j\big(\partial_{\underline x}f(\underline x)\partial_{\underline x}\big),\]
\begin{equation}\label{eq3}
\partial_{\underline x}\big(f(\underline x)e_j\big)\partial_{\underline x}=-2\partial_{x_j}\partial_{\underline x}f(\underline x)-\big(\partial_{\underline x}f(\underline x)\partial_{\underline x}\big)e_j,
\end{equation}
$j=1,\dots,m$.\hfill$\square$\vspace{0.24cm}

For a vector $\underline x$ and a $k$-vector $Y_k$, the inner and outer product between $\underline x$ and $Y_k$ are defined by (see \cite{DSS})
\[\underline x\bullet Y_k=\left\{\begin{array}{ll}\left[\underline xY_k\right]_{k-1}&\text{for}\;\;k\ge1\\0&\text{for}\;\;k=0\end{array}\right.\quad\quad\text{and}\quad\quad\underline x\wedge Y_k=\left[\underline xY_k\right]_{k+1}.\]
In a similar way $Y_k\bullet\underline x$ and $ Y_k\wedge\underline x$ are defined. We thus have that
\begin{alignat*}{1}
\underline xY_k&=\underline x\bullet Y_k+\underline x\wedge Y_k,\\
Y_k\underline x&=Y_k\bullet\underline x+Y_k\wedge\underline x,
\end{alignat*}
where also
\begin{alignat*}{1}
\underline x\bullet Y_k&=(-1)^{k-1}Y_k\bullet\underline x,\\
\underline x\wedge Y_k&=(-1)^kY_k\wedge\underline x.
\end{alignat*}
Let us now consider a $k$-vector valued function $F_k$ which is inframonogenic in the open set $\Omega\subset\mathbb R^m$. This is equivalent to say that $F_k$ satisfies in $\Omega$ the system
\[\left\{\begin{array}{ll}\partial_{\underline x}\bullet(\partial_{\underline x}\bullet F_k)&=0\\
\partial_{\underline x}\wedge(\partial_{\underline x}\bullet F_k)-\partial_{\underline x}\bullet(\partial_{\underline x}\wedge F_k)&=0\\
\partial_{\underline x}\wedge(\partial_{\underline x}\wedge F_k)&=0.
\end{array}\right.\]
In particular, for $m=2$ and $k=1$, a vector-valued function $\underline f=f_1e_1+f_2e_2$ is inframonogenic if and only if
\[\left\{\begin{array}{ll}
\partial_{x_1x_1}f_1-\partial_{x_2x_2}f_1+2\partial_{x_1x_2}f_2=0&\\
\partial_{x_1x_1}f_2-\partial_{x_2x_2}f_2-2\partial_{x_1x_2}f_1=0.&
\end{array}\right.\]
We now try to find particular solutions of the previous system of the form
\begin{align*}
f_1(x_1,x_2)&=\alpha(x_1)\cos(nx_2),\\
f_2(x_1,x_2)&=\beta(x_1)\sin(nx_2).
\end{align*}
It easily follows that $\alpha$ and $\beta$ must fulfill the system
\begin{align*}
\alpha^{\prime\prime}+n^2\alpha+2n\beta^\prime &=0\\
\beta^{\prime\prime}+n^2\beta+2n\alpha^\prime &=0.
\end{align*}
Solving this system, we get
\begin{align}
f_1(x_1,x_2)&=\big((c_1+c_2x_1)\exp(nx_1)+(c_3+c_4x_1)\exp(-nx_1)\big)\cos(nx_2),\label{sol1}\\
f_2(x_1,x_2)&=\big((c_3+c_4x_1)\exp(-nx_1)-(c_1+c_2x_1)\exp(nx_1)\big)\sin(nx_2)\label{sol2}.
\end{align}
Therefore, we can assert that the vector-valued function
\begin{multline*}
\underline f(x_1,x_2)=\big((c_1+c_2x_1)\exp(nx_1)+(c_3+c_4x_1)\exp(-nx_1)\big)\cos(nx_2)e_1\\
+\big((c_3+c_4x_1)\exp(-nx_1)-(c_1+c_2x_1)\exp(nx_1)\big)\sin(nx_2)e_2,\;\;c_j,n\in\mathbb R,
\end{multline*}
is inframonogenic in $\mathbb R^2$. Note that if $c_1=c_3$ and $c_2=c_4$, then
\begin{align*}
f_1(x_1,x_2)&=2(c_1+c_2x_1)\cosh(nx_1)\cos(nx_2),\\
f_2(x_1,x_2)&=-2(c_1+c_2x_1)\sinh(nx_1)\sin(nx_2).
\end{align*}
Since the functions (\ref{sol1}) and (\ref{sol2}) are harmonic in $\mathbb R^2$ if and only if $c_2=c_4=0$, we can also claim that not every inframonogenic function is harmonic.

Here is a simple technique for constructing inframonogenic functions from two-sided monogenic functions.

\begin{proposition}
Let $f(\underline x)$ be a two-sided monogenic function in $\Omega\subset\mathbb R^m$. Then $\underline xf(\underline x)$ and $f(\underline x)\underline x$ are inframonogenic functions in $\Omega$.
\end{proposition}
\textit{Proof.} It is easily seen that
\[\big(\underline xf(\underline x)\big)\partial_{\underline x}=\sum_{j=1}^m\partial_{x_j}\big(\underline xf(\underline x)\big)e_j=\underline x\big(f(\underline x)\partial_{\underline x}\big)+\sum_{j=1}^me_jf(\underline x)e_j=\sum_{j=1}^me_jf(\underline x)e_j.\]
We thus get
\[\partial_{\underline x}\big(\underline xf(\underline x)\big)\partial_{\underline x}=-\sum_{j=1}^me_j\big(\partial_{\underline x}f(\underline x)\big)e_j-2f(\underline x)\partial_{\underline x}=0.\]
In the same fashion we can prove that $f(\underline x)\underline x$ is inframonogenic.\hfill$\square$\vspace{0.24cm}

We must remark that the functions in the previous proposition are also harmonic. This may be proved using the following equalities
\begin{equation}\label{eq1}
\Delta_{\underline x}\big(\underline xf(\underline x)\big)=2\partial_{\underline x}f(\underline x)+\underline x\big(\Delta_{\underline x}f(\underline x)\big),
\end{equation}
\begin{equation}\label{eq2}
\Delta_{\underline x}\big(f(\underline x)\underline x\big)=2f(\underline x)\partial_{\underline x}+\big(\Delta_{\underline x}f(\underline x)\big)\underline x,
\end{equation}
and the fact that every monogenic function is harmonic. At this point it is important to notice that an $\mathbb R_{0,m}$-valued harmonic function is in general not inframonogenic. Take for instance $h(\underline x)e_j$, $h(\underline x)$ being an $\mathbb R$-valued harmonic function. If we assume that $h(\underline x)e_j$ is also inframonogenic, then from (\ref{eq3}) it may be concluded that $\partial_{\underline x}h(\underline x)$ does not depend on $x_j$. Clearly, this condition is not fulfilled for every harmonic function.

We can easily characterize the functions that are both harmonic and inframonogenic. Indeed, suppose that $h(\underline x)$ is a harmonic function in a star-like domain $\Omega\subset\mathbb R^m$. By the Almansi decomposition (see \cite{MR,R}), we have that $h(\underline x)$ admits a decomposition of the form
\[h(\underline x)=f_1(\underline x)+\underline x f_2(\underline x),\]
where $f_1(\underline x)$ and $f_2(\underline x)$ are left monogenic functions in $\Omega$. It is easy to check that
\[\partial_{\underline x}h(\underline x)=-mf_2(\underline x)-2\mathrm{E}_{\underline x}f_2(\underline x),\]
$\mathrm{E}_{\underline x}=\sum_{j=1}^mx_j\partial_{x_j}$ being the Euler operator. Thus $h(\underline x)$ is also inframonogenic in $\Omega$ if and only if $mf_2(\underline x)+2\mathrm{E}_{\underline x}f_2(\underline x)$ is right monogenic in $\Omega$. In particular, if $h(\underline x)$ is a harmonic and inframonogenic homogeneous polynomial of degree $k$, then $f_1(\underline x)$ is a left monogenic homogeneous polynomial of degree $k$ while $f_2(\underline x)$ is a two-sided monogenic homogeneous polynomial of degree $k-1$.

The following proposition provides alternative characterizations for the case of $k$-vector valued functions.

\begin{proposition}
Suppose that $F_k$ is a harmonic $($resp. inframonogenic$)$ $k$-vector valued function in $\Omega\subset\mathbb R^m$ such that $2k\ne m$. Then $F_k$ is also inframonogenic $($resp. harmonic$)$ if and only if one of the following assertions is satisfied:
\begin{itemize}
\item [{\rm(i)}] $F_k(\underline x)\underline x$ is left $3$-monogenic in $\Omega$;
\item [{\rm(ii)}] $\underline xF_k(\underline x)$ is right $3$-monogenic in $\Omega$;
\item [{\rm(iii)}] $\underline xF_k(\underline x)\underline x$ is biharmonic in $\Omega$.
\end{itemize}
\end{proposition}
\textit{Proof.} We first note that
\[e_je_Ae_j=\left\{\begin{array}{ll}(-1)^{\vert A\vert}e_A&\text{for}\quad j\in A,\\(-1)^{\vert A\vert+1}e_A&\text{for}\quad j\notin A,\end{array}\right.\]
which clearly yields $\sum_{j=1}^me_je_Ae_j=(-1)^{\vert A\vert}(2\vert A\vert-m)e_A$. It thus follows that for every $k$-vector valued function $F_k$,
\[\sum_{j=1}^me_jF_ke_j=(-1)^k(2k-m)F_k.\]
Using the previous equality together with (\ref{eq1}) and (\ref{eq2}), we obtain
\begin{align*}
\partial_{\underline x}\Delta_{\underline x}\big(F_k(\underline x)\underline x\big)&=2\partial_{\underline x}F_k(\underline x)\partial_{\underline x}+\big(\partial_{\underline x}\Delta_{\underline x}F_k(\underline x)\big)\underline x+(-1)^k(2k-m)\Delta_{\underline x}F_k,\\
\Delta_{\underline x}\big(\underline xF_k(\underline x)\big)\partial_{\underline x}&=2\partial_{\underline x}F_k(\underline x)\partial_{\underline x}+\underline x\big(\Delta_{\underline x}F_k(\underline x)\partial_{\underline x}\big)+(-1)^k(2k-m)\Delta_{\underline x}F_k,\\
\Delta_{\underline x}^2\big(\underline xF_k(\underline x)\underline x\big)&=4\Big(2\partial_{\underline x}F_k(\underline x)\partial_{\underline x}+(-1)^k(2k-m)\Delta_{\underline x}F_k+\big(\partial_{\underline x}\Delta_{\underline x}F_k(\underline x)\big)\underline x\\
&\quad\quad+\underline x\big(\Delta_{\underline x}F_k(\underline x)\partial_{\underline x}\big)\Big)+\underline x\big(\Delta_{\underline x}^2F_k(\underline x)\big)\underline x.
\end{align*}
The proof now follows easily.\hfill$\square$\vspace{0.24cm}

Before ending the section, we would like to make two remarks. First, note that if $m$ even, then a $m/2$-vector valued function $F_{m/2}(\underline x)$ is inframonogenic if and only if $F_{m/2}(\underline x)$ and $F_{m/2}(\underline x)\underline x$ are left $3$-monogenic, or equivalently, $F_{m/2}(\underline x)$ and $\underline xF_{m/2}(\underline x)$ are right $3$-monogenic. Finally, for $m$ odd the previous proposition remains valid for $\mathbb R_{0,m}$-valued functions.

\section{Fischer decomposition}

The classical Fischer decomposition provides a decomposition of arbitrary homogeneous polynomials in $\mathbb R^m$ in terms of harmonic homogeneous polynomials. In this section we will derive a similar decomposition but in terms of inframonogenic homogeneous polynomials. For other generalizations of the Fischer decomposition we refer the reader to \cite{CSV,BS,DSS,E,FK,MR,RM,S,SVN}.

Let $\mathsf{P}(k)$ ($k\in\mathbb N_0$) denote the set of all $\mathbb R_{0,m}$-valued homogeneous polynomials of degree $k$ in $\mathbb R^m$. It contains the important subspace $\mathsf{I}(k)$ consisting of all inframonogenic homogeneous polynomials of degree $k$.

An an inner product may be defined in $\mathsf{P}(k)$ by setting
\[\left\langle P_k(\underline x),Q_k(\underline x)\right\rangle_k=\left[\overline{P_k(\partial_{\underline x})}\,Q_k(\underline x)\right]_0,\quad P_k(\underline x),Q_k(\underline x)\in\mathsf{P}(k),\]
$\overline{P_k(\partial_{\underline x})}$ is the differential operator obtained by replacing in $P_k(\underline x)$ each variable $x_j$ by $\partial_{x_j}$ and taking conjugation.

From the obvious equalities
\begin{align*}
[\overline{e_ja}\,b]_0&=-[\overline ae_jb]_0,\\
[\overline{ae_j}\,b]_0&=-[\overline abe_j]_0,\quad a,b\in\mathbb R_{0,m},
\end{align*}
we easily obtain
\begin{align*}
\left\langle\underline xP_{k-1}(\underline x),Q_k(\underline x)\right\rangle_k&=-\left\langle P_{k-1}(\underline x),\partial_{\underline x}Q_k(\underline x)\right\rangle_{k-1},\\
\left\langle P_{k-1}(\underline x)\underline x,Q_k(\underline x)\right\rangle_k&=-\left\langle P_{k-1}(\underline x),Q_k(\underline x)\partial_{\underline x}\right\rangle_{k-1},
\end{align*}
with $P_{k-1}(\underline x)\in\mathsf{P}(k-1)$ and $Q_k(\underline x)\in\mathsf{P}(k)$. Hence for $P_{k-2}(\underline x)\in\mathsf{P}(k-2)$ and $Q_k(\underline x)\in\mathsf{P}(k)$, we deduce that
\begin{equation}\label{eq4}
\left\langle\underline xP_{k-2}(\underline x)\underline x,Q_k(\underline x)\right\rangle_k=\left\langle P_{k-2}(\underline x),\partial_{\underline x}Q_k(\underline x)\partial_{\underline x}\right\rangle_{k-2}.
\end{equation}

\begin{theorem}[Fischer decomposition]
For $k\ge2$ the following decomposition holds:
\[\mathsf{P}(k)=\mathsf{I}(k)\oplus\underline x\mathsf{P}(k-2)\underline x.\]
Moreover, the subspaces $\mathsf{I}(k)$ and $\underline x\mathsf{P}(k-2)\underline x$ are orthogonal w.r.t. the inner product $\left\langle\,,\right\rangle_k$.
\end{theorem}
\textit{Proof.} The proof of this theorem will be carried out in a similar way to that given in \cite{DSS} for the case of monogenic functions.

As $\mathsf{P}(k)=\underline x\mathsf{P}(k-2)\underline x\oplus\left(\underline x\mathsf{P}(k-2)\underline x\right)^\bot$ it is sufficient to show that \[\mathsf{I}(k)=\left(\underline x\mathsf{P}(k-2)\underline x\right)^\bot.\]
Take $P_k(\underline x)\in\left(\underline x\mathsf{P}(k-2)\underline x\right)^\bot$. Then for all $Q_{k-2}(\underline x)\in\mathsf{P}(k-2)$ it holds
\[\left\langle Q_{k-2}(\underline x),\partial_{\underline x}P_k(\underline x)\partial_{\underline x}\right\rangle_{k-2}=0,\]
where we have used (\ref{eq4}). In particular, for $ Q_{k-2}(\underline x)=\partial_{\underline x}P_k(\underline x)\partial_{\underline x}$ we get that $\partial_{\underline x}P_k(\underline x)\partial_{\underline x}=0$ or $P_k(\underline x)\in\mathsf{I}(k)$. Therefore $\left(\underline x\mathsf{P}(k-2)\underline x\right)^\bot\subset\mathsf{I}(k)$.

Conversely, let $P_k(\underline x)\in\mathsf{I}(k)$. Then for each $Q_{k-2}(\underline x)\in\mathsf{P}(k-2)$,
\[\left\langle\underline xQ_{k-2}(\underline x)\underline x,P_k(\underline x)\right\rangle_k=\left\langle Q_{k-2}(\underline x),\partial_{\underline x}P_k(\underline x)\partial_{\underline x}\right\rangle_{k-2}=0,\]
whence $P_k(\underline x)\in\left(\underline x\mathsf{P}(k-2)\underline x\right)^\bot$.\hfill$\square$\vspace{0.24cm}

By recursive application of the previous theorem we get:

\begin{corollary}[Complete Fischer decomposition]
If $k\ge2$, then
\[\mathsf{P}(k)=\bigoplus_{s=0}^{[k/2]}\underline x^s\mathsf{I}(k-2s)\underline x^s.\]
\end{corollary}

\subsection*{Acknowledgments}

D. Pe\~na Pe\~na was supported by a Post-Doctoral Grant of \emph{Funda\c{c}\~ao para a Ci\^encia e a Tecnologia}, Portugal (grant number: SFRH/BPD/45260/2008).

\end{document}